\documentclass{amsproc}
\usepackage{amssymb, euscript}
\newtheorem{theorem}{Theorem}[section]
\newtheorem{lemma}[theorem]{Lemma}
\newtheorem{corollary}[theorem]{Corollary}
\theoremstyle{definition}
\newtheorem{definition}[theorem]{Definition}
\newtheorem{example}[theorem]{Example}

\theoremstyle{remark}
\newtheorem{remark}[theorem]{Remark}

\numberwithin{equation}{section}

\def\Cal{\mathcal}

\def\H{{\Cal H}}

\def\P{{\Cal P}}

\def\F{{\Cal F}}

\def\I{{\Cal I}}

\def\S{\EuScript{S}}

\def\tr{{\hbox{\rm tr}}}

\def\Ma{\frM_{n,m}}

\def\Mmm{\frM_{m,m}}

\def \lv{{\bf \lam}}
\def \nv{{\bf n}_0}
\def \mn{{\bf m}_0}
\def \kv{{\bf k}_0}

\def\Z{\mathcal{Z}}

\def\f0{f_0}
\def\Fc0{\varphi_0}

\def\I_k {I_{-}^{k/2}}
\def\I+k {I_{+}^{k/2}}

\def\vnm{V_{n,m}}

\def\gnm{G_{n,m}}

\def\bbr{{\Bbb R}}

\def\bbc{{\Bbb C}}

\def\rank{{\hbox{\rm rank}}}

\def\const{{\hbox{\rm const}}}

\def\det{{\hbox{\rm det}}}

\def\part{\partial}
\def\intl{\int\limits}
\def\b{\beta}

\def\Gam{\Gamma}
\def\Om{\Omega}
\def\a{\alpha}
\def\om{\omega}

\def\Del{\Delta}

\def\vp{\varphi}

\def\g{\gamma}

\def\Lam{\Lambda}
\def\sig{\sigma}
\def\lam{\lambda}
\def\lv{{\boldsymbol \lam}}
\def\mv{{\boldsymbol \mu}}
\def\av{{\boldsymbol \a}}
\def\bv{{\boldsymbol \b}}
\def\z{\zeta}

\def\t{\tau}

\font\frak=eufm10

\def\fr#1{\hbox{\frak #1}}

\def\frL{\fr{L}}        
\def\frM{\fr{M}}

\def\const{{\hbox{\rm const}}}

\def\det{{\hbox{\rm det}}}

\def\p{\Om}
\def\gm{\Gamma_m}

\def\part{\partial}
\def\intl{\int\limits}
\def\b{\beta}

\def\Gam{\Gamma}
\def\Om{\Omega}
\def\a{\alpha}
\def\cpm{\overline\p}

\newcommand{\be}{\begin{equation}}
\newcommand{\ee}{\end{equation}}
\newcommand{\bea}{\begin{eqnarray}}
\newcommand{\eea}{\end{eqnarray}}
\newcommand{\Bea}{\begin{eqnarray*}}
\newcommand{\Eea}{\end{eqnarray*}}

\begin{document}

\title[The Composite Cosine Transform] {The Composite Cosine Transform on
 the Stiefel Manifold and Generalized Zeta Integrals}

\author{ E. Ournycheva}
 \address{ Department  of Mathematics, Bar Ilan  University,
 52900 Ramat Gan,   ISRAEL} \email{ournyce@macs.biu.ac.il}

\author{ B. Rubin}
 \address{Department of Mathematics, Louisiana State
University,  Baton Rouge, LA, 70803, USA, \newline
 Institute of Mathematics, Hebrew University, Jerusalem 91904,
 ISRAEL}
\email{borisr@math.lsu.edu}

\thanks{ The work  was supported in part
by the Edmund Landau Center for Research in Mathematical Analysis
and Related Areas, sponsored by the Minerva Foundation (Germany),
and the Abraham and Sarah Gelbart Research Institute for
Mathematical Sciences. }

\date{January 2, 2005}



\begin{abstract} The
$\lam$-cosine transform on the unit sphere is defined by
\[(T^{\lam} f)(u)=\int_{S^{n-1}} f(v) |v \cdot u|^{\lam}
\, dv, \qquad u \in S^{n-1},\] and has many applications. We
introduce a new integral transform $T^{\lv} f, \lv \in \bbc^m$,
which generalizes the previous one for functions on the Stiefel and
 Grassmann manifolds. We call it the composite cosine transform, by
taking into account that its kernel agrees with the composite power
function of the cone of positive definite symmetric matrices. Our
aim is to describe  the set of all $\lv \in \bbc^m$ for which
$T^{\lv}$ is injective on the space of integrable functions. We
obtain the precise description of this set  in some important
 cases, in particular, for $\lam$-cosine transforms on
Grassmann manifolds. The main tools are the classical Fourier
analysis of functions of matrix argument and the relevant zeta
integrals.
\end{abstract}

\maketitle

\section{Introduction}
\setcounter{equation}{0}

The classical cosine transform (also known as the Blaschke-Levy
transform) on the unit sphere $S^{n-1}$ in $\bbr^n$ is defined by
\be\label{t11}(T f)(u)=\intl_{S^{n-1}} f(v) |v \cdot u| \, dv,
\qquad u \in S^{n-1},\ee where  $f$ is an integrable even function
on $S^{n-1}$, and $ v \cdot u$ is the usual inner product. This
transform and its generalization \be\label{t1}(T^{\lam}
f)(u)=\intl_{S^{n-1}} f(v) |v \cdot u|^{\lam} \, dv,\ee arise in
diverse areas of mathematics, in particular, in PDE, the Fourier
analysis, integral geometry,  and the Banach space theory; see
\cite{Es}, \cite{Ga},
 \cite{GH1}, \cite{GG},
\cite {Ko},  \cite{Ru2}-\cite{Ru4}, \cite{Sa}, \cite{Schn},
\cite{Se}. Operators (\ref{t1}) have been
  investigated in detail using the following two approaches.
The first one employs the Fourier transform technique \cite {Se},
\cite {Ko}, \cite {Ru2}, and relies
 on the formula \be\label{ff} \intl_{\bbr^{n}} \frac{f\big({x \over
|x|}\big )}{|x|^{n+\lam}} e^{ix \cdot y} dx = c_{\lam,n}
|y|^{\lam} (T^{\lam}f) \big({y \over |y|}\big ), \quad
c_{\lam,n}=\const, \ee which should be properly interpreted. The
second approach is based on decomposition in spherical harmonics.
Namely, by the Funk-Hecke formula, \be\label{sh}
T^{\lam}P_k=c\,\mu_k (\lam)\, P_k,\ee \be\label{mu1}
c=2\pi^{(n-1)/2}(-1)^{k/2},\qquad \mu_k (\lam)=
 \frac{\Gam\Big(\frac{\lam+1}{2}\Big)\,
\Gam\Big(\frac{k-\lam}{2}\Big)}{\Gam\Big(-\frac{\lam}{2}\Big) \,
\Gam\Big(\frac{\lam+k+n}{2}\Big)}, \ee for each  homogeneous
harmonic polynomial $P_k (x)$ of even degree $k$ restricted to the
unit sphere; see \cite{Ru2}-\cite{Ru4}, \cite{Sa}. The
Fourier-Laplace multiplier $\mu_k (\lam)$ provides complete
information about properties of $T^{\lam}$.

In the last two decades a considerable attention  was attracted to
generalizations of  $T$ and $T^{\lam}$ for functions on the
Grassmann manifold $G_{n,m}$  of $m$-dimensional linear subspaces of
$\bbr^n$. We recall, that if  $\eta\in G_{n,m}$, $ \xi \in G_{n,l},
\; l \ge m$, and $[\eta |\xi ]$ is  the $m$-dimensional volume of
the parallelepiped spanned by the orthogonal projection of a generic
orthonormal coordinate frame in $\eta$ onto $\xi$, then, by
definition, \be\label{ccgr} (T^{\lam} f)(\xi)=\intl_{G_{n,m}}
f(\eta) \, [\eta | \xi]^{\lam} \, d\eta\ee (we use the same notation
as in (\ref{t1})). For $l>m$, the operator (\ref{ccgr}) represents
the composition of the similar one over $G_{n,l}$ and the
corresponding Radon transform acting from $G_{n,m}$ to
 $G_{n,l}$ (see, e.g., \cite {A}, \cite{GR}). Thus injectivity results
 for $T^\lam$ in this case can be easily derived from those for the Radon transform
 (see \cite{GR} and references therein) and the case $l=m$.
 Owing to this, in the following we
assume $l=m$. This case bears the basic features of the operator
family (\ref{ccgr}).

The investigation of operators (\ref{ccgr}) for $\lam=1$ was
initiated in stochastic geometry (processes of flats) by Matheron
\cite{Mat1}, \cite{Mat2} who conjectured that $T^1$ is injective
as well as its prototype (\ref{t11}). Matheron's conjecture was
disproved by Goodey and Howard \cite{GH1} who used the idea of
Gluck and Warner \cite{GW} to interpret  the Grassmann manifold
$G_{4,2}$ as the direct product $S^2 \times S^2$ of 2-spheres. For
higher dimensions, the result then follows by induction. Operators
$T^\lam$ for $\lam = 0,1,2,\dots$ were studied  in  \cite[p.
117]{GH2}, where, by using reduction to $G_{4,2}$, it was  proved
that $T^\lam$ is non-injective for such $\lam$; see also
\cite{GHR}, \cite{Sp1}, \cite{Sp2}, \cite{Gr}. The range of the
$\lam$-cosine transform was studied by Alesker and Bernstein
\cite{AB} for $\lam=1$ and by Alesker \cite{A} for all complex
$\lam$, who invoked deep results from the representation theory.

In the present article we develop a new approach to operators
$T^{\lam}$. Our argument differs essentially from that in the cited
papers. For technical reasons, we prefer
 to deal with $O(m)$ right-invariant functions on the
 Stiefel manifold $\vnm$ of orthonormal $m$-frames
 rather than with functions on  $G_{n,m}$. This
 does not change the essence of the matter and leads to the
 following equivalent definition:
\be\label{tlll}(T^{\lam} f)(u)=\intl_{\vnm} f(v) |\det(v'u)|^{\lam}
\, dv,\qquad u \in \vnm, \ee  where ``$\;{}'\;$" stands for the
transposed matrix, and the product $v' u$ is understood in the sense
of matrix mulilication.
 Then we
 regard (\ref{tlll}) as a member of the more general analytic family
\be \label{tnff}(T^{\lv} f)(u)=\intl_{\vnm} f(v) \, (u'v v'u)^{\lv}
\, dv, \qquad u\in\vnm, \ee
 where $
\lv=(\lam_1,\ldots, \lam_m) \in \bbc^m$ and $(\cdot )^\lv$ denotes
the composite power function of the cone of positive definite $m
\times m$  matrices; see Section 2.2. We call $T^{\lv}f$ {\it the
composite cosine transform} of a function $f$ on $\vnm$. The general
intention is to obtain a higher-rank analog of the formula
(\ref{sh}), evaluate the corresponding multiplier $\mu_k (\lv)$
explicitly in terms of gamma functions, and use it for examination
of $T^{\lv}$. The particular case $\lv=(1, \ldots , 1)$ corresponds
to Matheron's operator. We do not realize this project in full
generality here and leave this work for future publications.
However, our approach enables us to obtain the precise description
of those $\lv$ for which $T^{\lv}$ is injective in the following
important cases (a) $2m \le n, \; $$ \lv=(\lam_1,\ldots, \lam_m)\in
\bbc^m$, and (b) $\lam_1=\dots =\lam_m=\lam \in \bbc$, provided that
$T^{\lv}f$ and $T^{\lam}f$ exist in the usual Lebesgue sense as
absolutely convergent integrals.

The essence of our approach is that we
  apply the classical Fourier transform technique  to obtain
 higher-rank analogs of  (\ref{ff}) and (\ref{sh}). In the
 first case we assume $f$ to be an  arbitrary integrable
  function on $\vnm$, and in the second one $P_k$ stands for the
 restriction to $\vnm$ of the corresponding  $O(m)$ right-invariant determinantally
homogeneous harmonic polynomial  on the space of  $n \times m$
matrices.  Different aspects of harmonic analysis based on
implementation of such polynomials were studied in \cite{Herz},
\cite{Str}, \cite{TT}, and in a series of publications related to
 group representations. In the present article we
obtain a higher-rank
 copy of (\ref{sh}) with the multiplier $\mu_k (\lv)$ explicitly
 expressed in terms of the gamma functions associated with the cone
 of positive definite $m\times m$ matrices.  In the particular case
 $\lam_1=\dots =\lam_m$
  the main result reads as follows.
\begin{theorem}  Let $n> m \ge2, \; f \in L^1(\gnm)$. Then
$(T^{\lam} f)(\xi)$ is finite for almost all $\xi \in \gnm$ if and
only if $Re \, \lam > -1$. For such $\lam$, the operator $T^{\lam}$
is injective on $L^1(\gnm)$ if and only if $\lam \neq
0,1,2,\dots$.\end{theorem} In particular, we show that if  $\lam$ is
a non-negative integer and $2m<n$, then $T^\lam$, having been
written in the form (\ref{tlll}), annihilates  all $O(m)$
right-invariant determinantally homogeneous harmonic polynomials
$P_k(x)$  of degree $k> Re \,\lam +m-1$.

 A by-product of our investigation is  a
 functional equation for the generalized zeta
integrals with additional ``angle component" $f(v), \; v \in \vnm$;
see (\ref{zeta}), (\ref{zeta*}). This result is of independent
interest. The main references related to zeta integrals can be found
in \cite{B}, \cite{BSZ}, \cite{FK}, \cite{Kh2}, \cite{Ru5}. The
equation obtained below is, in fact,
  far-reaching  higher-rank modifications of
(\ref{ff}) in the language of Schwartz distributions. In cited
papers  they were obtained for $f\equiv 1$. The case $f=P_k$ was
considered  in \cite{Cl} in the context of Jordan algebras. The
proof presented below  is much simpler than that in \cite{Cl}
(adapted to our case) and employs the idea from \cite{Kh2} to derive
the result for $\lam \in \bbc$ from the more general one for $\lv
\in \bbc^m$.

The paper is organized as follows. Section 2 contains the necessary
background material. In Section 3 we introduce the generalized zeta
integrals. Section 4 plays the central role in the article and
establishes connection between zeta integrals and composite cosine
transforms. The results of Section 4 are applied in Section 5 to
study injectivity of the composite cosine transform.

{\bf Acknowledgement.} We are thankful to Dr. S.P. Khekalo for
sharing with us his results \cite{Kh2}.

\section{Preliminaries}

\setcounter{equation}{0}

We establish our notation and recall basic facts that will be used
throughout the paper. The main references are \cite{FK},
\cite{Gi}, \cite{OR}, \cite{T}.
\subsection{Notation}
  Let $\frM_{n,m}$ be the
space of real matrices $x=(x_{i,j})$ having $n$ rows and $m$
 columns. We identify $\frM_{n,m}$
 with the real Euclidean space $\bbr^{nm}$ and set $dx=\prod^{n}_{i=1}\prod^{m}_{j=1}
 dx_{i,j}$ for the Lebesgue measure on $\Ma$. If $n \ge m$, then
 $\frM_{n,m}^0$ stands for the set of all matrices $x \in \frM_{n,m}$
 of rank $m$. This set has a full measure in $\frM_{n,m}$.
 In the following,
  $x'$ denotes the transpose of  $x$, $I_m$ is the identity $m \times m$
  matrix,  $0$ stands for zero entries. Given a square matrix $a$,  we denote by $|a|$ the absolute value of
 the determinant of $a$, and  by $\tr (a)$  the trace of $a$, respectively.

 Let  $\p=\P_m$ be the cone
of positive definite symmetric matrices $r=(r_{i,j})_{m\times m}$
with the elementary volume $ dr=\prod_{i \le j} dr_{i,j}$, and let
$\cpm$ be  the closure of $\p$, that is the set of  all positive
semi-definite $m\times m$ matrices. For $r\in\p$ ($r\in\cpm$) we
write $r>0$ ($r\geq 0$). Given $s_1$ and  $s_2$ in $\cpm$, the
inequality $s_1 > s_2$  means $s_1 - s_2 \in
 \p$.
  If $a\in\cpm$ and $b\in\p$, then   $\int_a^b f(s) ds$ denotes
 the integral over the set
$$
 \{s : s \in \p, \, a<s<b \}=\{s : s-a \in \p, \, b-s \in \p \}.
$$
The group $G=GL(m,\bbr)$ of
 real non-singular $m \times m$ matrices $g$ acts  on
 $\p$   by the rule $r \to \tau_g (r)=grg'$ so that $\tau_{g_1} \tau_{g_2}
  =\tau_{g_1 g_2}; \; g_1, g_2 \in G$.  The corresponding $G$-invariant
 measure on $\p$ is  \be\label{2.1}
  d_{*} r = |r|^{-(m+1)/2} dr, \qquad |r|=\det (r),  \ee \cite[p.
  18]{T}. The group $G$ is transitive on
 $\p$ but not simply transitive. The transitivity retains if we restrict
 to the subgroup $T_m$ of upper triangular matrices with positive diagonal elements.
     This subgroup is simply transitive.  The subgroup of lower
     triangular  matrices with positive diagonal elements also acts simply transitively  on
 $\p$. Each $r \in
\p$ has a unique representation $r=t't, \; t=(t_{i,j}) \in T_m$,
so that \be \label{ya-tr}dr=2^m \prod_{j=1}^m t_{j,j}^{m-j+1}
dt_{j,j}\, dt_\ast, \quad dt_{*}=\prod_{i<j}
 dt_{i,j};\ee \cite[p.
  39]{T}. An alternative representation  reads
 $r=\tau\tau'$, $\tau \in T_m$. To connect both
 representations, let
\be\label{om} r_\ast =\om
r\om,\quad \om=\left[\begin{array}{ccccc} 0 & {}   & {}    & 1 \\
                              {} & {}  & {.}    & {} \\
                              {} & {.}   & {}   & {} \\
                               1 & {}   & {}    & 0

\end{array} \right],\quad \om^2=I_m. \ee If $r \in
\p$ and $r_\ast=t't, \; t \in T_m$, then $r=\tau\tau'$, where
$\tau=\om t'\om\in T_m$.

We use a standard notation $O(n)$   and $SO(n)$  for the group of
real orthogonal $n\times n$ matrices and its connected component
of the identity, respectively. The invariant measure on  $SO(n)$
is normalized to be of total mass 1.
 The Schwartz space  $\S=\S(\Ma)$ is
 identified with  the respective space on $\bbr^{nm}$.

 The Fourier transform  of a
function $f\in L^1(\Ma)$ is defined by \be\label{ft} (\F
f)(y)=\intl_{\Ma} e^{{\rm tr(iy'x)}} f (x) dx,\qquad y\in\Ma \;
.\ee This is the usual Fourier transform on $\bbr^{nm}$, and  the
relevant Parseval equality  reads \be\label{pars} (\F f, \F
\vp)=(2\pi)^{nm} \, (f,  \vp),\ee where
$$(f, \vp)=\intl_{\Ma} f(x) \overline{\vp(x)} \, dx.$$

\begin{lemma}\label{12.2} {\rm (see, e.g.,
 \cite[pp. 57--59]{Mu} )}.\hskip10truecm

\noindent
 {\rm (i)} \ If $ \; x=ayb$, where $y\in\Ma, \; a\in  GL(n,\bbr)$, and $ b \in  GL(m,\bbr)$, then
 $dx=|a|^m |b|^ndy.$\\
 {\rm (ii)} \ If $ \; r=qsq'$, where $s\in \P_m$ and $ q\in  GL(m,\bbr)$,
  then $dr=|q|^{m+1}ds.$ \\
  {\rm (iii)} \ If $ \; r=s^{-1}$,  $s\in \p$,   then $r\in
  \p$
  and $dr=|s|^{-m-1}ds.$
\end{lemma}

\subsection{The composite power function}
Given $r=(r_{i,j})\in\p$, let $\Del_0(r)=1$, $\Del_1(r)=r_{1,1}$,
$\Del_2(r)$, $\ldots$, $\Del_m(r)=|r|$ be the corresponding
principal minors which are strictly positive \cite[p. 586]{Mu}.
For $\lv=(\lam_1,\dots,\lam_m)\in\bbc^m$, the composite power
function of the cone  $\p$  is defined by
 \bea\label{pf} r^{\lv}&=&\prod\limits_{i=1}^m \left[\frac{\Del_i (r)}{\Del_{i-1} (r)}\right]^{\lam_i/2}\\
&=& \Del_1 (r)^{\frac{\lam_1 -\lam_2}{2}}  \ldots \Del_{m-1}
(r)^{\frac{\lam_{m-1} -\lam_m}{2}} \Del_m (r)^{\frac{\lam_m}{2}}.
\eea
\begin{remark}
 In the case $m=1$, the function (\ref{pf}) becomes
$r^{\lam/2}$. This notational confusion is easily resolved if one
takes into account the equality (\ref{r-tr}) below.
\end{remark}

The composite power  functions associated to homogeneous cones
were introduced by S. Gindikin \cite{Gi}. We also refer to
\cite{FK}, \cite{Kh2},  and \cite{T}, where (\ref{pf}) is written
in different notation.

We denote \[|\lv|=\lam_1+\cdots+\lam_m.\]
 If $\lam_1/2,\dots,\lam_m/2$ are integers and
$\lam_1\geq\lam_2\geq\cdots\geq\lam_m$, then $r^{\lv}$ is a
polynomial of degree $|\lv|/2$. In the special case $
\lam_1=\ldots=\lam_m=\lam$ we add the subscript $0$ so that
$$ \lv_0= (\lam, \ldots, \lam) \; (\in \bbc^m)
$$
and
 \be\label{det}r^{{\lv_0}}=|r|^{\lam/2}.\ee If $r=t't, \; t \in
T_m$, then $\Del_i(r)=[\Del_i(t)]^2=\prod_{j=1}^{i}t_{j,j}^2\;$,
and therefore,  \be\label{r-tr}
 r^{{\lv}}=\prod\limits_{j=1}^m t_{j,j}^{\lam_j}\equiv \pi_\lv(t).
\ee The expression $\pi_\lv(t)$ is a multiplicative character of
 $T_m$, so that  \be
\pi_\lv(t_1t_2)=\pi_\lv(t_1) \, \pi_\lv(t_2), \qquad
\pi_\lv(t^{-1})=\pi_{-\lv}(t). \ee

For $\lv=(\lam_1,\dots,\lam_m)$, let
$\lv_\ast=(\lam_m,\dots,\lam_1)$ be the reverse vector;
$(\lam_\ast)_j=\lam_{m-j+1}$. To each $r\in\p$ we associate the
matrix $r_\ast=\om r\om $ (see (\ref{om}))  with the components
$(r_\ast)_{i,j}=r_{m-j+1,\, m-i+1}$.

\begin{lemma} \label {uh} Let $\lv, \mv\in\bbc^m$;  $r\in\p$.  Then  \bea\label{pr1}
r^{{\lv+\mv}}&=&r^{{\lv}}\;r^{{\mv}}, \quad
r^{{\lv+\av_0}}=r^{\lv}|r|^{\a/2}, \quad \av_0=(\a,\dots,\a);
\\ \label{pr6} \; (t'rt)^{{\lv}}&=&(t't)^{\lv}\;r^{{\lv}},\quad t\in T_m;
\\\label{pr4}
 r^{{\lv_\ast}}&=&(r^{-1})_\ast^{-\lv},\quad
(r^{-1})^{{\lv}}=r_\ast^{{-\lv_\ast}};\\
\label{pr5} \qquad (cr)^{\lv}&=&c^{|\lv|/2}r^{\lv} ,\quad
(cr)^{\lv}_*=c^{|\lv|/2}r_*^{\lv} , \quad c>0. \eea
\end{lemma}
\begin{proof}
These statements (up to notation) may be found in different
sources \cite{Gi}, \cite{FK}, \cite{Kh2}, \cite{T}. For the sake
of completeness, we outline the proof. The property  (\ref{pr1})
is clear in view of (\ref{r-tr}) and (\ref{det}). To prove
(\ref{pr6}), it suffices to set $r=\tau'\tau$, $\tau\in T_m$, and
make use of (\ref{r-tr}). To prove (\ref{pr4}), let $r=t't$, $t\in
T_m$. Then by (\ref{r-tr}),
$$ r^{{\lv_\ast}}=\prod\limits_{j=1}^m t_{j,j}^{\lam_{m-j+1}}.$$
On the other hand, since  $r^{-1}=t^{-1} (t')^{-1}$, and
$(r^{-1})_\ast=\om r^{-1}\om=(\om t^{-1}\om )(\om t^{-1}\om )',$
by taking into account the equalities
$$
t^{-1}=\left[\begin{array}{ccccc} t_{1,1}^{-1}    & {}  & {}  & {} \\
                              {} & {.}    &{ }  & {*} \\
                              {0} & {}   & {.}   & {} \\
                              {} & {}    & {}  & t_{m,m}^{-1}

\end{array} \right],\qquad \om t^{-1}\om =\left[\begin{array}{ccccc} t_{m,m}^{-1}    & {}  & {}  & {} \\
                              {} & {.}    &{ }  & {0} \\
                              {*} & {}   & {.}   & {} \\
                               {} & {}    & {}  & t_{1,1}^{-1}

\end{array} \right],
$$
\cite[p. 580 (vii)]{Mu}, we have \bea\nonumber
(r^{-1})_\ast^{-\lv}&=&\left (\left[\begin{array}{ccccc} t_{m,m}^{-1}    & {}  & {}  & {} \\
                              {} & {.}    &{ }  & {0} \\
                              {*} & {}   & {.}   & {} \\
                              {} & {}    & {}  & t_{1,1}^{-1}

\end{array} \right]\left[\begin{array}{ccccc} t_{m,m}^{-1}    & {}  & {}  & {} \\
                              {} & {.}    &{ }  & {*} \\
                              {0} & {}   & {.}   & {} \\
                               {} & {}    & {}  & t_{1,1}^{-1}

\end{array} \right]\right)^{-\lv}\\\nonumber
&=&\prod\limits_{j=1}^m t_{m-j+1,\,
m-j+1}^{\lam_{j}}=\prod\limits_{j=1}^m t_{j,j}^{\lam_{m-j+1}}.\eea
This gives $ r^{{\lv_\ast}}=(r^{-1})_\ast^{-\lv}$. Replacing $r$
by $r^{-1}$ and $\lv_\ast$ by $\lv$, we obtain the second equality
in (\ref{pr4}). The first equality in (\ref{pr5}) follows
immediately from (\ref{pf}). The second equality is a consequence
of the first one and (\ref{pr4}):
\[ (cr)^{\lv}_*= (c^{-1} r^{-1})^{-\lv_*} = c^{|\lv|/2}(r^{-1})^{-\lv_*}=
c^{|\lv|/2} r_*^{\lv}.\]
\end{proof}

\subsection{Gamma  function of the cone $\p$}

 The  gamma function of the cone $\p$ is defined
by \be\label{gf} \Gam_{\Omega}  (\lv) =\intl_{\Omega} r^{\lv}
e^{-{\rm tr} (r)} d_{*} r.
 \ee
This integral converges absolutely
 if and only if $Re \, \lam_j>j-1$, and represents a product of
 ordinary $\Gamma$-functions:
\be\label{sv}\Gam_{\Omega}
(\lv)=\pi^{m(m-1)/4}\prod\limits_{j=1}^{m} \Gam ((\lam_j-
j+1)/2),\ee see, e.g., \cite[p. 123]{FK}.
\begin{lemma}
For $s\in\p$ and $Re \, \lam_j>j-1$, \be\label{eq1} \intl_{\Omega}
r^{\lv} e^{-{\rm tr} (rs)} d_{*} r= \Gam_{\Omega} (\lv) \,
s_\ast^{{-\lv_\ast}}.\ee
\end{lemma}
\begin{proof}
This equality  is known (\cite[p.  23 ]{Gi}, \cite[p. 124 ]{FK},
\cite{Kh2}). By (\ref{pr4}), it is equivalent to \be
\label{eq1.1}\intl_{\Omega} r^{\lv} e^{-{\rm tr} (rs^{-1})} d_{*}
r= \Gam_{\Omega} (\lv) \, s^{\lv}.\ee Assuming $s=t't$, $t\in
T_m$, and changing variable $r=t'\rho t$, for the left-hand side
of (\ref{eq1.1}) we obtain
$$
\intl_{\Omega} (t'\rho t)^{\lv} e^{-{\rm tr} (\rho)} d_{*} \rho
\stackrel{\rm (\ref{pr6})}{=}\Gam_{\Omega} (\lv)
(t't)^{\lv}=\Gam_{\Omega} (\lv) s^{\lv}.
$$
\end{proof}

An important particular case of (\ref{gf}) is  the Siegel integral
\be\label{sgf} \Gam_{m}  (\lam) =\intl_{\Omega} |r|^{\lam} e^{-{\rm
tr} (r)} d_{*} r=\pi^{m(m-1)/4}\prod\limits_{j=0}^{m-1} \Gam (\lam -
j/2),\ee which converges absolutely
 if and only if $Re \,\lam > (m-1)/2$. By (\ref{det}),
\be\label{sv1} \Gam_{\Omega}  (\lv_0)=\gm(\lam/2),\quad
\lv_0=(\lam,\dots,\lam). \ee

\subsection{Stiefel  manifolds}

For $n\geq m$, let $\vnm= \{v \in \frM_{n,m}: v'v=I_m \}$
 be the Stiefel manifold of orthonormal $m$-frames in $\bbr^n$.
 We fix the invariant measure $dv$ on
 $\vnm$ [Mu, p. 70]
  normalized by \be\label{2.16} \sigma_{n,m}
 \equiv \intl_{\vnm} dv = \frac {2^m \pi^{nm/2}} {\gm
 (n/2)}\; \ee and denote  $d_\ast v=\sig^{-1}_{n,m} dv$.
The polar decomposition on $\Ma$ is defined according to the
following lemma; see, e.g.,  \cite[pp. 66, 591]{Mu}, \cite{Ma}.
\begin{lemma}\label{l2.3} If $x \in \frM_{n,m}, \; \rank (x)=m, \; n \ge m$,
then \be \label{pol} x=vr^{1/2}, \qquad v \in \vnm,   \qquad r=x'x
\in\p,\ee and $dx=2^{-m} |r|^{(n-m-1)/2} dr dv$.
\end{lemma}

 A modification of Lemma \ref{l2.3} in terms of upper
triangular matrices  reads as follows.
\begin{lemma}\label{sph} {\rm (\cite{P}, \cite{Ru5})} If
$x \in \frM_{n,m}, \; \rank (x)=m,\;  n \ge m$, then
  \[
x=ut, \qquad u \in \vnm,   \qquad t \in T_m,\] and
$$
dx=\prod\limits_{j=1}^m t_{j,j}^{n-j} dt_{j,j}dt_*dv, \qquad
dt_*=\prod\limits_{i<j} dt_{i,j}.
$$
\end{lemma}

\section{The generalized zeta integrals}

By taking into account the polar decomposition $x=vr^{1/2}$  (see
Lemma \ref{l2.3}), we introduce the following generalized zeta
integrals (or zeta distributions): \bea\label{zeta}
\Z(\phi,\lv,f)&=&\intl_{\Ma}
r^{\lv} \, f(v) \, \overline{\phi(x)} \,  dx=(r^\lv f, \phi), \\
\label{zeta*}\Z_*(\phi,\lv,f)&=&\intl_{\Ma} r_*^{\lv} \, f(v) \,
\overline{\phi(x)}  \, dx=(r_*^\lv f, \phi). \eea Here $r=x'x$,
$v=xr^{-1/2}$, $f$ is a fixed integrable function on $\vnm$, and
$\phi \in \S(\Ma)$ is a test function.

The following particular cases are worth  mentioning.

{\bf $1^0$}.  $\lv=\lv_0=(\lam, \ldots, \lam)$.

In this case  \be\label{lr} \Z(\phi,\lv, f)=\Z_*(\phi,\lv,
f)=\intl_{\Ma} (\det \, r)^{\lam/2}  \, f(v)\, \overline{\phi(x)}
\, dx.\ee For $f\equiv 1$, zeta integrals of this type were
studied in \cite {Ge}, \cite {Kh1}; see also
 \cite {Ru5} and references therein.

{\bf $2^0$}. $m=1$.

In this case $r=x'x=|x|^2$ and (\ref{lr}) reads \be
\Z(\phi,\lam,f)=\intl_{\bbr^n} |x|^\lam f \big (\frac{x}{|x|}\big
) \, \overline{ \phi(x)} \, dx. \ee Distributions of this form are
well known in analysis and integral geometry; see, e.g.,
\cite{Se}, \cite{Es},  \cite{Ko}, \cite{Ru2}.

{\bf $3^0$}. $\lv \in \bbc^m, \; f\equiv 1$.

This case was explored in \cite{FK} (in the context of Jordan
algebras) by invoking the relevant $K$-Bessel functions; see also
\cite{Ra}, \cite{Kh2}, \cite{B}, \cite{BSZ}, and references therein.

{\bf $4^0$}. $\lv \in \bbc^m, \; f$ is a determinantlally homogeneous harmonic
 polynomial (see Section 4.2).

 In this case, integrals (\ref{zeta}) and (\ref{zeta*})   were studied
 in \cite {Cl} using the argument close to \cite{FK}.

Most of the publications mentioned above were focused on evaluation
of the Fourier transform of the corresponding zeta distributions. This transform is
 realized in the form of the relevant functional equation of the Parseval type.

Our nearest goal is to investigate convergence  of integrals
(\ref{zeta}) and (\ref{zeta*}), and their analyticity in the
$\lv$-variable. We denote  \bea\label{abs} {\bf
\Lam}&=&\{\lv\in\bbc^m :
Re\,\lam_j>j-n-1 \quad \forall j=1,\dots, m \}, \\
\label{pset}
{\bf \Lam}_0&=&\{\lv\in\bbc^m : \lam_j=j-n-l \quad \text{\rm for some}\\
 &{}&j \in \{1,\dots, m\}, \quad \text{\rm and } \quad l \in \{1,3,5,\dots \} \}.
 \nonumber
 \eea

\begin{lemma}\label{l3.1}  The integrals (\ref{zeta}) and
(\ref{zeta*}) are  absolutely convergent if and only if $\lv \in
{\bf\Lam}$ and extend
 as  meromorphic functions of $ \lv$
with the  polar set ${\bf \Lam}_0 $ .  The normalized zeta
integrals \be\label{nzi} \Z^0(\phi,\lv,f)=\frac{
\Z(\phi,\lv,f)}{\Gam_\Om(\lv+\nv)}, \qquad
\Z^0_*(\phi,\lv,f)=\frac{\Z_*(\phi,\lv,f)}{\Gam_\Om(\lv+ \nv)},
\ee $\nv=(n, \ldots, n)$, are   entire functions of $\lv$.
 \end{lemma}

\begin{proof} Let us consider the integral (\ref{zeta}). We set
 $ x=vt, \quad v \in
\vnm$,   $ t \in T_m$, and make use of Lemma \ref{sph}. By taking
into account that $x'x=t't$ and $t(t't)^{-1/2} \in O(m)$, owing to
(\ref{r-tr}), we obtain \be\label{dis1} \Z(\phi,\lv,f)=
\intl_{\bbr^m_+} F(t_{1,1},\dots , t_{m,m}) \prod\limits_{j=1}^m
t_{j,j}^{\lam_j+n-j} dt_{j,j}\, , \ee where
\[
F(t_{1,1},\dots , t_{m,m}) = \intl_{\bbr^{m(m-1)/2}}
dt_*\intl_{\vnm} f(v)\overline{\phi(v t)}\,d v, \quad
dt_*=\prod\limits_{i<j} dt_{i,j}.\] Since $F$ extends as an even
Schwartz function in each argument, it can be written as
$$
F(t_{1,1},\dots , t_{m,m}) =F_0(t^2_{1,1},\dots , t^2_{m,m}),
$$
where $F_0\in\S(\bbr^m)$ (use, e.g., Lemma 5.4 from \cite[p.
56]{Tr}). Replacing $t_{j,j}^2$ by $s_{j,j}$, we represent
(\ref{dis1}) as a direct product of one-dimensional distributions
\be\label{z-h} \Z(\phi,\lv,f)=(\prod\limits_{j=1}^m
(s_{j,j})_+^{(\lam_j+n-j-1)/2},\;  F_0(s_{1,1},\dots , s_{m,m})).\ee
It follows that  the integral (\ref{zeta}) is absolutely convergent
provided $Re\,\lam_j>j-n-1$, i.e., $\lv \in {\bf\Lam}$. The
condition  $\lv \in {\bf\Lam}$ is strict. This claim becomes clear
if we choose $f\equiv 1$ and $\phi(x)=e^{-{\rm tr} (x'x)}$ which
give \be \label{ex}\Z(\phi,\lv,f)=\intl_{\Ma} (x'x)^{\lv} \,
e^{-{\rm tr} (x'x)}\,dx=2^{-m}\sig_{n,m}\Gam_\Om(\lv+\nv).\ee
Furthermore, since $(s_{j,j})_+^{(\lam_j+n-j-1)/2}$ extends as a
meromorphic distribution with the only poles  $\lam_j= j-n-1, j-n-3,
\dots \;$, then, by the fundamental Hartogs theorem \cite {Sh}, the
function $ \lv \to \Z(\phi,\lv,f)$ extends
 as a meromorphic function
with the  polar set ${\bf \Lam}_0$. By the same reason, a direct
product of the normalized distributions $
(s_{j,j})_+^{(\lam_j+n-j-1)/2}/\Gam ((\lam_j+n-j+1)/2)$ is an
entire function of $\lv$.

Let us consider the integral (\ref{zeta*}). By changing variable
$x=y\om$ where $\om$ is the matrix (\ref{om}), we obtain $
r_*=(x'x)_*=(\om y'y\om)_*=y'y$. Hence (set $y=u\t, \; u \in \vnm,
\; \t \in T_m$),  \bea \Z_*(\phi,\lv,f) &=&\intl_{\Ma} (y'y)^{\lv}
\, f(y\om (\om
y' y \om)^{-1/2}) \, \overline{\phi(y\om)} \, dy \nonumber \\
&=& \intl_{\bbr^m_+} \Phi(\t_{1,1},\ldots , \t_{m,m})
\prod\limits_{j=1}^m \t_{j,j}^{\lam_{j}+n-j} d\t_{j,j}\,
\nonumber\eea where, as above (note that $\t\om \, (\om \t'\t
\om)^{-1/2} \in O(m)$), \bea \Phi(\t_{1,1},\ldots , \t_{m,m}) &=&
\intl_{\bbr^{m(m-1)/2}} d\t_*\intl_{\vnm} f(u) \, \overline{\phi(u
\t\om)}\,d u\nonumber \\
&=&\Phi_0(\t^2_{1,1},\dots , \t^2_{m,m}), \qquad
\Phi_0\in\S(\bbr^m). \nonumber \eea This gives
\[\Z_*(\phi,\lv,f)=(\prod\limits_{j=1}^m
(s_{j,j})_+^{(\lam_{j}+n-j-1)/2},\; \Phi_0(s_{1,1},\dots ,
s_{m,m})),\] and the result follows as in the previous case.
\end{proof}

\section{Connection between zeta integrals and the composite cosine transform}

\subsection{The composite cosine transform}
Let $\lv=(\lam_1,\ldots, \lam_m)\in\bbc^m$,  $f$ be an integrable
function on $\vnm$. We  introduce the following family of
intertwining operators \be \label{tnf}(T^{\lv} f)(u)=\intl_{\vnm}
f(v) \, (u'v v'u)^{\lv} \, dv, \qquad u\in\vnm,\qquad n>m, \ee
which commute with the left action of the orthogonal group $O(n)$.
We call $(T^{\lv} f)(u)$ {\it the
 composite cosine transform  of } $f$. In the particular case
 $\lam_1=\ldots=\lam_m=\lam$, the
integral (\ref{tnf})  represents the usual $\lam$-cosine transform
 \be\label{tll}(T^{\lam}
f)(u)=\intl_{\vnm} f(v) |\det( v' u)|^{\lam} \, dv.\ee

The function $(T^{\lv} f)(u)$ extends to all matrices $y \in \Ma$ of
rank $m$ if we set $y=ut$, $ u \in \vnm$,   $ t \in T_m$. Indeed, by
(\ref{pr6}),  \be \label{tf1}(T^{\lv} f)(y)=r^{\lv}(T^{\lv} f)(u)\ee
where $r^{\lv}=(t't)^{\lv}=(y'y)^{\lv}$ is ``the radial part" of
$(T^{\lv} f)(y)$.

\begin{theorem}\label {exi}  For $f \in L^1(\vnm)$, the integral $(T^{\lv} f)(u)$
converges absolutely for almost all $u \in \vnm$ if and only if $Re
\, \lam_j > j-m-1$ for all $j=1,2, \ldots ,m$ and represents an
analytic function of  $\lv$ in this domain. For such $\lv$, the
linear operator $T^{\lv} $ is bounded on $L^1(\vnm)$. \end{theorem}

This statement follows immediately from the following lemma which is
of independent interest.

\begin{lemma}\label {exl} Let $u,v \in \vnm$. Then

 \be \label {ave}I \equiv \intl_{\vnm} (u'v v'u)^{\lv} \,
du=\intl_{\vnm} (u'v v'u)^{\lv} \, dv=\frac{2^m \,
\pi^{nm/2}}{\Gam_m (m/2)}\, \frac{\Gam_\Om(\lv
+\mn)}{\Gam_\Om(\lv+\nv)}\,. \ee This integral converges absolutely
if and only if $Re \, \lam_j > j-m-1$ for all $j=1,2, \ldots ,m$.
\end{lemma}
\begin{proof} The first equality in (\ref{ave}) becomes clear
if we write both integrals as those over the  group $SO(n)$ (see
Lemma 2.4 in \cite{GR}). These integrals  are, in fact, constant
with respect to the corresponding exterior variables.  Thus, one
can write \be \label {ax} I=\intl_{\vnm} (u'v_0 v_0'u)^{\lv} \,
du, \qquad v_0=
\left[\begin{array} {c}  I_{m} \\
0 \end{array} \right] \in \vnm. \ee To evaluate $I$ we introduce
 an auxiliary integral
\be \label {aax} A=\intl_{\Ma} (x'v_0 v_0'x)^{\lv} \, e^{-{\rm tr}
(x'x)} \, dx,\ee
 and transform it in two different ways. Let first $$x=\left[\begin{array} {c}
  a \\ b \end{array} \right], \qquad a \in \frM_{m,m}, \quad b\in
  \frM_{n-m,m}.$$ Then $v_0'x=a, \, x'x= a'a +b'b$, and we have
\be \label {q1} A=A_1 A_2, \quad A_1=\intl_{\frM_{m,m}} (a'a)^{\lv}
\, e^{-{\rm tr} (a'a)} \, da, \quad A_2= \intl_{\frM_{n-m,m}}
e^{-{\rm tr} (b'b)} \, db.\ee By Lemma \ref{l2.3}, (\ref{pr1}),
(\ref{gf}), and (\ref{2.16}), we obtain \be \label {q2} A_1=2^{-m}
\sig_{m,m} \intl_{\Omega} r^{\lv+\mn} e^{-{\rm tr} (r)} d_{*} r=
\frac{ \pi^{m^2/2}\, \Gam_\Om(\lv +\mn)}{\Gam_m (m/2)}\ee provided
$Re \, \lam_j
> j-m-1$,  $j=1,2, \ldots ,m$. The last condition is sharp and provides
the ``only if" part in the lemma and in the Theorem \ref{exi}. For
$A_2$ we have
\[ A_2=\left
( \, \intl_{-\infty}^{\infty}e^{-s^2}\ ds \right )^{m(n-m)}=
\pi^{m(n-m)/2}.\] Thus \be \label {q22} A=\frac{ \pi^{nm/2}\,
\Gam_\Om(\lv +\mn)}{\Gam_m (m/2)}, \qquad Re \, \lam_j
> j-m-1.\ee

On the other hand, by setting $x=ut, \, u \in \vnm,   \, t \in T_m$,
owing to Lemma \ref{sph}, we obtain
\[A=\intl_{T_m} e^{-{\rm tr} (t't)} \, d \mu (t)\intl_{\vnm}
(t'u'v_0 v_0'ut)^{\lv} \, du,\]
\[ d \mu (t)=\prod\limits_{j=1}^m t_{j,j}^{n-j} \, dt_{j,j} \, dt_*, \qquad
dt_*=\prod\limits_{i<j} dt_{i,j}.\] By (\ref{pr6}), one can write
\be \label {q3}A=BI,\ee where $I$ is our integral (\ref{ax})  and
\bea B&=&\intl_{T_m} (t't)^\lv e^{-{\rm tr} (t't)}
\,d\mu(t)=\nonumber
\\ &=&\prod\limits_{j=1}^m \, \intl_0^\infty t_{j,j}^{\lam_j+n-j}\,
e^{-t_{j,j}^2}\, dt_{j,j} \times \prod\limits_{i<j}\,
\intl_{-\infty}^\infty e^{-t_{i,j}^2}\, dt_{i,j} \nonumber \\
&=&2^{-m} \,\pi ^{m(m-1)/4} \, \prod\limits_{j=1}^m \Gam \Big
(\frac{\lam_j+n-j+1}{2}\Big )\nonumber \\ &=&2^{-m}
\,\Gam_\Om(\lv+\nv), \qquad Re \, \lam_j > j-n-1.\nonumber \eea
Combining this with (\ref{q22}) and (\ref {q3}), we obtain
\[ 2^{-m} \,\Gam_\Om(\lv+\nv) \, I=\frac{
\pi^{nm/2}\, \Gam_\Om(\lv +\mn)}{\Gam_m (m/2)}\] provided that $Re
\, \lam_j
> j-m-1$ for all $j=1,2, \ldots ,m$. This gives (\ref{ave}).
\end{proof}
\begin{corollary}\label{le} Let
$\lam_1=\ldots=\lam_m=\lam, \quad Re\; \lam > -1$. Then
\be\label{mnv}
 \intl_{V_{n, m}} |\det( v' u)|^{\lam} \,
dv=\dfrac{2^m\, \pi^{nm/2}\, \Gam_m
\big(\frac{m+\lam}{2}\big)}{\Gam_m \big(\frac{m}{2}\big) \, \Gam_m
\big(\frac{n+\lam}{2}\big)}. \ee  \end{corollary}

If we set $\lam =1$ in (\ref{mnv}) and replace $dv$ by the
normalized measure $d_*v$ then the resulting formula bears an
 important geometrical meaning. Namely, let $\xi$ and $\eta$ be
 $ m$-dimensional linear subspaces of $\bbr^n$,
 i.e., $ \xi, \eta \in G_{n,m}$, and let
  $u$ and $v$ be   coordinate frames
of    $\xi$ and $\eta$, respectively. Then $|\det( v' u)|=[\xi |
\eta]$ is  the  volume of the projection onto $\eta$ of the
parallelepiped spanned by $u$. The corresponding averaged volume is
valuated as \be\label{clam1} \intl_{V_{n, m}} |\det( v' u)| \,
d_*v=\intl_{G_{n,m}}  [\xi | \eta] \, d\eta= \dfrac{\Gam_m
\big(\frac{n}{2}\big) \, \Gam_m \big(\frac{m+1}{2}\big)}{\Gam_m
\big(\frac{m}{2}\big) \, \Gam_m \big(\frac{n+1}{2}\big)}.\ee

 \subsection{The basic functional equation}
 Below we establish connection between zeta integrals and composite cosine
 transforms. This can be done in the form of a
  functional equation which is, in fact, the usual  Parseval
 equality in the framework of the corresponding Fourier analysis.    We start
with the following
\begin{lemma}
 Let $f $ be an integrable $O(m)$  right-invariant  function on
 $\vnm$; $x=vr^{1/2}, \; \psi (x)=f(v) |r|^{(m-n)/2}$. For $s\in\p$ and
$\phi\in\S(\Ma)$, \bea\label{eq-psi} &{}&\frac{
|s|^{-m/2}}{\sig_{m,m}}\intl_{\Ma}\overline{(\F \phi)(y)}
dy\intl_{\vnm}f(v) \, e^{-{\rm tr}(\pi v' ys^{-1}y'v)} d
v\\\nonumber &{}& \quad = \, (2\pi)^{(n-m)m}\intl_{\Ma} \psi (x)
\, e^{-{\rm tr}(xsx'/4\pi)} \, \overline{\phi(x)} \, dx.\eea
\end{lemma}
\begin{proof} By the  Parseval equality (\ref{pars}), it suffices
to find the Fourier transform of the function $ \psi_s(x)=\psi
(x)\, e^{-{\rm tr}(xsx'/4\pi)}$. In  polar coordinates we have
\bea (\F\psi_s)(y)&=&2^{-m}\intl_{\vnm}f(v)\, d
v\intl_{\p}|r|^{-1/2} e^{{\rm
tr}(iy'vr^{1/2}-r^{1/2}sr^{1/2}/4\pi)}\;dr \nonumber \\
&{}&  \text{\rm (replace $v$ by $v\g$, $ \g\in O(m)$, and
integrate in $\g$)}\nonumber
\\&=&\!\!2^{-m}\!\!\intl_{\vnm}\!\!f(v)\, d v\intl_{O(m)} \!\!d\g
\intl_{\p}|r|^{-1/2} e^{{\rm tr}(iy'v\g r^{1/2}-
r^{1/2}sr^{1/2}/4\pi)}\;dr\nonumber
\\&=&\!\!\frac{1}{\sig_{m,m}}\intl_{\vnm}\!\!f(v)\, d
v \intl_{\Mmm} e^{{\rm tr}(iy'vz-zsz'/4\pi)}\;dz. \nonumber \eea
The inner integral can be evaluated by the  formula \be
\intl_{\Mmm} e^{{\rm tr}(i \z'z)}\,  e^{-{\rm tr}(zsz'/4\pi)}
dz=(2\pi)^{m^2}|s|^{-m/2}e^{-{\rm tr}(\pi \z s^{-1}\z')}, \ee
where
 $\z=v'y$,   \cite[p. 481]{Herz}. This gives
\[(\F\psi_s)(y)=\frac{(2\pi)^{m^2}
|s|^{-m/2}}{\sig_{m,m}}\intl_{\vnm}f(v)e^{-{\rm tr}(\pi v'
ys^{-1}y'v)} d v, \] and (\ref{eq-psi}) follows.
\end{proof}

Let us compute the Fourier transform of the function
\be\label{fil} \vp_{\lv}(x)= \frac{ r_*^{-\lv_*-\nv}}{
\Gam_\Om(-\lv_*)} \, f(v)\ee where $x=vr^{1/2}, \; \nv=(n,\ldots,
n)$.
\begin{theorem}\label{gen}
Let $f $ be  an integrable $O(m)$  right-invariant function on
 $\vnm$, $\lv\in\bbc^m$.
Then \be\label{eq5}
(\F\vp_{\lv})(y)=\frac{c_{\lv}}{\Gam_\Om(\lv+\mn)} \, (T^{\lv}
f)(y), \qquad c_{\lv}=2^{-|\lv|}\pi^{m^2/2}/\sig_{m,m},\ee
 in the sense of $\S'$-distributions, i.e.
\bea \label{eq9}\frac{c_{\lv}}{\Gam_\Om(\lv+\mn)} \, ( T^{\lv}
f,\; \F\phi)&=&(2\pi)^{nm} \, (\vp_\lv, \; \phi)\\ &\equiv&
(2\pi)^{nm} \, \Z^0_*(\phi,-\lv_*-\nv, f) \nonumber \eea for each
$\phi\in \S(\Ma)$.
\end{theorem}
\begin{proof}
The equality (\ref{eq9}) demonstrates an intimate interrelation
between zeta integrals and cosine transforms. To prove (\ref{eq9}),
we multiply (\ref{eq-psi}) by $s^{\lv +\mn}$ and integrate against
$d_\ast s$. We obtain \bea\label{eq2} &{}&\frac{
1}{\sig_{m,m}}\intl_{\Ma}\overline{(\F \phi)(y)}\,
dy\intl_{\vnm}f(v)\, A_1(v'y) \, d v \\&{}& \qquad = \;
(2\pi)^{(n-m)m}\intl_{\Ma} \psi (x) \,  A_2(x) \, \overline{\phi(x)}
\, dx,\nonumber \eea where $\psi (x)=f(v) |r|^{(m-n)/2}$,
\bea\nonumber A_1(\z)&=&\intl_{\p} s^{\lv +\mn}|s|^{-m/2} e^{-{\rm
tr}(\pi \z s^{-1}\z')}d_\ast s,\qquad \z\in\Mmm,\\\nonumber
A_2(x)&=&\intl_{\p}s^{\lv+\mn} e^{-{\rm tr}(xsx'/4\pi)}d_\ast
s,\qquad x\in\Ma. \eea Let us compute $A_1$ and $A_2$. By
(\ref{pr1}), setting $a=\pi \z' \z$, $r=s^{-1}$, we obtain
 \bea A_1(a)&=&\intl_{\p} (r^{-1})^{ \lv}
e^{-{\rm tr}(a r )}d_\ast r \nonumber \\ &\stackrel{\rm
(\ref{pr4})}{=}&\intl_{\p} r_\ast^{\;-\lv_*} e^{-{\rm tr}(a r
)}d_\ast r \qquad \text{\rm ($r_\ast =\om r\om$; see (\ref{om}))}
\nonumber \\
&=&\intl_{\p} r^{\;-\lv_*} e^{-{\rm tr}(a \om r\om)}d_\ast r
\nonumber \\
&=&\intl_{\p} r^{\;-\lv_*} e^{-{\rm tr}(a_* r )}d_\ast r
\stackrel{\rm (\ref{eq1})}{=}\Gam_{\Omega} (-\lv_*) \, a^{\lv}
\eea provided $Re \,(\lam_*)_j=Re \, \lam_{m-j+1}<1-j$, or
$Re\,\lam_j<j-m$. Finally, by  (\ref{pr5}), we have
$$
A_1(\z)= \pi^{|\lv|/2} \Gam_{\Omega} (-\lv_*) (\z'
\z)^{\lv},\qquad Re\,\lam_j<j-m.
$$
Furthermore, by (\ref{eq1}) and (\ref{pr5}),
$$A_2(x)=(4\pi)^{(|\lv|+m^2)/2}\Gam_{\Omega} (\lv+\mn) (x' x)_*^{-\lv_* -\mn}, \qquad Re\,\lam_j>j-m-1.$$

Hence, (\ref{eq2}) reads \bea\label{eq7} &{}&\frac{ \pi^{|\lv|/2}
}{\sig_{m,m}\Gam_{\Omega} (\lv+\mn)}\intl_{\Ma}\overline{(\F
 \phi)(y)} \, (T^{\lv}
f)(y) \, dy \\ &{}& \qquad =
\frac{(2\pi)^{(n-m)m}(4\pi)^{(|\lv|+m^2)/2}}{\Gam_{\Omega}
(-\lv_*)} \intl_{\Ma} r_*^{-\lv_* -\nv} \, f(v)
\,\overline{\phi(x)} \, dx,\nonumber \\ &{}&\label{eq77} \qquad =
(2\pi)^{(n-m)m}(4\pi)^{(|\lv|+m^2)/2}\, \Z^0_*(\phi,-\lv_*-\nv,
f),  \eea \be\label{dom1}j-m-1<Re\,\lam_j<j-m.\ee By Lemma
\ref{l3.1}, the expression (\ref{eq77}) extends as an entire
function of $\lv$. Hence,   it can be regarded as analytic
continuation of the integral  (\ref{eq7}) outside of the domain
(\ref{dom1}), and the result follows.
\end{proof}

\begin{example} Let $m\ge 1, \; \lam_1=\ldots=\lam_m=\lam, \; |x|_m
=\det(x'x)^{1/2}$. Then
\[\vp_{\lv}(x)= \frac{|x|_m^{-\lam-n}}{ \Gam_m (-\lam/2)} \,
f(x(x'x)^{-1/2}), \qquad (T^{\lam} f)(y)=\intl_{\vnm} f(v) |\det(v'
y)|^{\lam} \, dv, \] and we have \be\label{rnm}
(\F\vp_{\lv})(y)=\frac{2^{-\lam m} \, \pi^{m^2/2}
}{\sig_{m,m}\Gam_m((\lam +m)/2)} \, (T^{\lam} f)(y).\ee
\end{example}
\begin{example} Let $m=1$. Then for $x, y \in \bbr^n \setminus
\{0\}$,
\[\vp_{\lv}(x)= \frac{|x|^{-\lam-n}}{ \Gam (-\lam/2)} \,
f\Big(\frac{x}{|x|}\Big), \qquad (T^{\lam} f)(y)=\intl_{S^{n-1}}
f(v) |v \cdot y|^{\lam} \, dv. \]  In this case \be\label{rn1}
(\F\vp_{\lv})(y)=\frac{2^{-1-\lam}\pi^{1/2}}{\Gam((\lam +1)/2 )}
\, (T^{\lam} f)(y).\ee
\end{example}

\subsection{The case of  homogeneous
polynomials}

Let $P_k(x)$ be a polynomial  on $\Ma$ which is harmonic (as a
function on $\bbr^{nm}$) and determinantally homogeneous of degree
$k$, so that \be P_k(xg)=\det (g)^k P_k(x),\quad \forall g\in
GL(m,\bbr),\quad x\in\Ma. \ee The latter means that $P_k$ is a usual
homogeneous harmonic polynomial of degree $km$ on $\bbr^{nm}$.
Theorem \ref{gen} can be essentially strengthened  if we choose $f$
to be the restriction of $P_k (x)$ onto $\vnm$. According to
(\ref{fil}), for $x=vr^{1/2}$, we denote \be\label{fik} \vp_{\lv,
k}(x)= \frac{ r_*^{-\lv_*-\nv}}{ \Gam_\Om(-\lv_*)} \, P_k(v).\ee

\begin{theorem} Let $y=vr^{1/2}, \; v \in \vnm, \; r \in \p$.
For all $\lv \in \bbc^m$,  \be\label{eq50}
(\F\vp_{\lv,k})(y)=\frac{d_{\lv}\,
\Gam_\Om(\kv-\lv_*)}{\Gam_\Om(-\lv_*) \, \Gam_\Om(\lv+\kv+\nv)} \,
P_k(v)\, r^{\lv},\ee \[ d_{\lv}=2^{-|\lv|}\pi^{nm/2}i^{km}, \] in
the sense of $\S'$-distributions, i.e., \bea \frac{d_{\lv}\,
\Gam_\Om(\kv-\lv_*)}{\Gam_\Om(\lv+\kv+\nv)} \,( P_k(v)\,
r^{\lv},\; \F\phi)&=&(2\pi)^{nm}(\vp_{\lv,k}, \; \phi) \nonumber
\\&=& (2\pi)^{nm}\,\Z^0_*(\phi,-\lv_*-\nv, P_k)\nonumber \eea
for each $\phi\in \S(\Ma)$.\end{theorem}
\begin{proof}
By the Hecke identity \cite[p. 87]{St}, \be \intl_{\Ma}P_k(x)
e^{-{\rm tr}(\pi x'  x)}e^{{\rm tr}(2\pi i
y'x)}\;dx=i^{km}P_k(y)e^{-{\rm tr}(\pi y'  y)},\ee or (replace $x$
by $xs^{1/2}/2\pi$ and $y$ by $ys^{-1/2}, \; s \in \Om$) $$
\intl_{\Ma}P_k(x) e^{-{\rm tr}( xs x'/4\pi)} \, e^{{\rm tr}(i
y'x)}\;dx=(2\pi)^{(k+n)m}i^{km}|s|^{-k-n/2} P_k(y)e^{-{\rm tr}(\pi
 y s^{-1}y')}.$$ By the Parseval equality, \bea
\label{eq-2}&&|s|^{-k-n/2}\intl_{\Ma} P_k(y)e^{-{\rm tr}(\pi y
s^{-1}y')}\, \overline{\F\phi(y)}\, dy\\\nonumber &=&(2\pi
i)^{-km}\intl_{\Ma} P_k(x) e^{-{\rm tr}( xs x'/4\pi)}\,
\overline{\phi(x)}\, dx. \eea We multiply  (\ref{eq-2}) by
$s^{\lv+\nv+\kv}$ and integrate against $d_\ast s$. This gives  \[
\intl_{\Ma} P_k(y)I_1(y)\, \overline{\F\phi(y)}\, dy=(2\pi
i)^{-km}\intl_{\Ma} P_k(x) I_2(x)\, \overline{\phi(x)}\, dx, \]
 where
\[I_1(y)=\intl_{\p} s^{\lv-\kv} e^{-{\rm tr}(\pi y
s^{-1}y')}d_\ast s,\quad I_2(x)=\intl_{\p}s^{\lv+\nv +\kv}
e^{-{\rm tr}( xs x'/4\pi)}d_\ast s. \] As in the proof of Theorem
\ref{gen}, we have
  $$ I_1(y)=\pi
^{(|\lv|-km)/2}\Gam_{\Omega} (\kv-\lv_*) (y' y)^{ \lv-\kv}, \quad
Re\,\lam_j <j+k-m,
$$  and $$
I_2(x)=(4\pi)^{(|\lv|+nm+km)/2}\Gam_{\Omega} (\lv+\nv+\kv) (x'
x)_*^{-\lv_* -\nv-\kv},$$ $$Re\,\lam_j >j-n-k-1.$$ Hence, if $
j-n-k-1<Re\,\lam_j<j+k-m$, then \bea \label{eq8} &&\Gam_{\Omega}
(\kv-\lv_*) \intl_{\Ma} P_k(y)\, (y' y)^{ \lv-\kv}
\,\overline{\F\phi(y)}\, dy\\\nonumber &=&c_\lv\Gam_{\Omega}
(\lv+\kv+\nv) \intl_{\Ma} P_k(x)\, (x' x)_*^{-\lv_* -\nv-\kv}
\,\overline{\phi(x)}\, dx, \eea $
c_\lv=i^{-km}2^{nm+|\lv|}\pi^{nm/2}. $ By taking into account
homogeneity of $P_k$ (i.e., $P_k(vr^{1/2})= P_k(v)\,|r|^{k/2}$),
and setting $\psi_{\lv,k} (x)=P_k(v)r^{\lv}, \;
 x=vr^{1/2}$, this can be written as \bea
\frac{d_{\lv}\, \Gam_\Om(\kv-\lv_*)}{\Gam_\Om(-\lv_*) \,
\Gam_\Om(\lv+\kv+\nv)} \, (\psi_{\lv,k},\F\phi)&=&(2\pi)^{nm}\,
 (\vp_{\lv,k},\phi) \nonumber\\
 &=&\label{eq555}(2\pi)^{nm} \,\Z^0_*(\phi,-\lv_*-\nv, P_k), \eea
 \be\label{dom} j-n-k-1<Re\,\lam_j<j+k-m,
\ee
 $d_{\lv}=2^{-|\lv|}\pi^{nm/2}i^{km}$.
Since the domain(\ref{dom}) is not void for all $k=0,1, 2,
\ldots$, and the normalized zeta integral (\ref{eq555}) is  an
entire function of $\lv$ (use Lemma \ref{l3.1}),  the result
follows.
\end{proof}

\begin{example} Let $m\ge 1, \; \lam_1=\ldots=\lam_m=\lam$,
$$\vp_{\lam, k}(x)=\frac{|x|_m^{-\lam-n}}{ \Gam_m (-\lam/2)} \, P_k
(x(x'x)^{-1/2}).$$ Then \be \label{ft2} (\F\vp_{\lam,k})(y)=\frac{
d_{\lam} \,  \Gam_m((k\!-\!\lam)/2)}{\Gam_m(-\lam/2) \,
\Gam_m((\lam\!+\!k\!+\!n)/2)} \, |y|_m^{\lam} \, P_k
(y(y'y)^{-1/2}),\ee
\[ d_{\lam}=2^{-\lam m}\pi^{nm/2}i^{km}.\]
\end{example}
\begin{example} Let $m=1$,  $$\vp_{\lam, k}(x)=\frac{|x|^{-\lam-n}}{ \Gam (-\lam/2)} \,
P_k \Big(\frac{x}{|x|}\Big).$$ Then \be \label{ft1}
(\F\vp_{\lam,k})(y)\!=\!\frac{d_{\lam} \,
\Gam((k-\lam)/2)}{\Gam(-\lam/2) \, \Gam((\lam+k+n)/2)} |y|^{\lam}
P_k \Big(\frac{y}{|y|}\Big),\ee \[
d_{\lam}\!=\!2^{-\lam}\pi^{n/2}i^{k}.\]
\end{example}

  \section{Injectivity of the composite cosine transform}
We denote by $\frL$ the set of all $\lv=(\lam_1, \dots ,\lam_m) \in
\bbc^m$ satisfying $Re \, \lam_j >j-m-1$ for all $j=1, \dots , m$.
By Theorem \ref{exi}, $(T^{\lv} f)(v)$ exists a.e. on $\vnm$ for
every $f \in L^1(\vnm)$ if and only if $\lv \in \frL$. In the
following we focus only on this domain, although other $\lv$ can
also be treated using analytic continuation. We also restrict our
consideration to the space $L^\flat (\vnm)$ which consists of $O(m)$
right-invariant integrable functions on $\vnm$. This space is
isomorphic to the space $L^1(G_{n,m})$ of integrable functions on
the Grassmann manifold $G_{n,m}$.

\begin{theorem}\label {main} Let $n> m$ and $\lv=(\lam_1, \dots \lam_m) \in \frL$.
 If \be \label {uhh}\lam_j +m-j \neq 0, 2, 4,\dots
\quad \text{\rm for all} \quad j=1, \dots , m,\ee then the operator
$$(T^{\lv} f)(u)=\intl_{\vnm} f(v) \, (u'v v'u)^{\lv} \, dv$$
is injective on $L^\flat (\vnm)$. If $2m\le n$ and (\ref{uhh}) fails
then $T^{\lv}$ is non-injective.
\end{theorem}

To prove this theorem, minor preparation is still needed.
\begin{definition} Following \cite{Herz}, we call a polynomial $P_k(x)$
on $\Ma$ an $H$-polynomial of degree $k$ if it is $O(m)$
right-invariant, harmonic, and determinantally homogeneous of degree
$k$. We denote by $\H_k$ the space of all such polynomials.
\end{definition}

\begin{lemma} \label{inff} Let $P_k \in \H_k$,
\be\label{mul} \mu_k(\lv)=\frac{\Gam_\Om(\lv+\mn)\,
\Gam_\Om(\kv-\lv_*)}{\Gam_\Om(-\lv_*) \, \Gam_\Om(\lv+\kv+\nv)},
\qquad \lv\in \bbc^m.\ee If $\lv$
 does not belong to the polar set of $\Gam_\Om(\lv+\mn)$, then
\be\label{cmp} (T^{\lv} P_k)(vr^{1/2})=c \,\mu_k(\lv)\,
P_k(v)r^{\lv}, \qquad c=\pi^{m(n-m)/2} \, i^{km}\,\sig_{m,m},\ee
in the sense of $\S'$-distributions.
\end{lemma}
\begin{proof} Let us  compare (\ref{eq5}) and (\ref{eq50}), assuming  $f$ to be the
restriction of $P_k$ onto $\vnm$.
 For all $\lv \in
\bbc^m$, we obtain \be\frac{c_{\lv}}{\Gam_\Om(\lv+\mn)} \, (T^{\lv}
P_k)(vr^{1/2})
 =\frac{d_{\lv}\,
\Gam_\Om(\kv-\lv_*)}{\Gam_\Om(-\lv_*) \, \Gam_\Om(\lv+\kv+\nv)} \,
P_k(v)\, r^{\lv},\ee
\[ c_{\lv}=2^{-|\lv|}\pi^{m^2/2}/\sig_{m,m}, \qquad
d_{\lv}=2^{-|\lv|}\pi^{nm/2}i^{km},\] in the sense of
$\S'$-distributions. If we exclude all $\lv$ belonging to the
polar set of $\Gam_\Om(\lv+\mn)$, then we get (\ref{cmp}).
\end{proof}
\begin{corollary} \label{cr}If $\lv \in \frL$ and $P_k \in \H_k$,
then for all $v \in \vnm$,\be\label{lmu}(T^{\lam} P_k)(v)=c \,\mu_k
(\lv)\, P_k(v),\ee  $c$ and $\mu_k (\lv)$ having the same meaning as
in (\ref{cmp}).
\end{corollary}
\begin{proof}  The function
\[\Gam_\Om(\lv+\mn)=\pi^{m(m-1)/4}\,\prod\limits_{j=1}^{m}
\Gam \Big(\frac{\lam_j+m- j+1}{2}\Big)\] has no poles in the domain
$ \frL$, and therefore, by (\ref{cmp}),\be\label{for10} ((T^{\lv}
P_k)(vr^{1/2}), \phi)=c \,\mu_k(\lv)\, (P_k(v)r^{\lv}, \phi) \ee for
all $\phi (y)\equiv \phi (vr^{1/2}) \in \S(\Ma)$. Now we choose
$\phi (y)=\chi (r)\, \psi (v)$ where $\chi (r)$ is a non-negative
$C^\infty$ cut-off function supported away from the boundary of $\p$
and $\psi (v)$ is a $C^\infty$
 function on $\vnm$. By passing to polar coordinates,
  from (\ref{for10}) and (\ref{tf1}) we obtain
 $$c_\chi \,\intl_{\vnm} [(T^{\lv}
P_k)(v)-c \,\mu_k(\lv)\,P_k (v)]\, \psi (v)\, dv =0, \qquad
c_\chi=\const \neq 0.$$ This implies (\ref{lmu}).
\end{proof}

\begin{remark} A simple computation shows that for $k=0$, the
equality (\ref{cr}) coincides with (\ref{ave}). Note also that if
$m=1$ then $P_k \in \H_k$ is necessarily of even degree, and we have
\be\label{comp2} (T^{\lam} P_k)(v)=2\pi^{(n-1)/2}(-1)^{k/2}\,
 \frac{\Gam\Big(\frac{\lam+1}{2}\Big)\,
\Gam\Big(\frac{k-\lam}{2}\Big)}{\Gam\Big(-\frac{\lam}{2}\Big) \,
\Gam\Big(\frac{\lam+k+n}{2}\Big)} \,  P_k (v).\ee This coincides
with the known formula (\ref{sh}) for spherical harmonics.
\end{remark}

\begin{remark} An important question is, do there exist $H$-polynomials of a given
 degree $k$? Note that for $n=m$ we have exactly two such
polynomials, namely, $P_0(x)\equiv 1$ and $P_1(x)=\det (x)$. The
following statement is due to Herz \cite{Herz}, p. 484: {\it For
$2m\le n$ there exist $H$-polynomials of every degree $k$}. It is
also known \cite{TT}, p. 27, that for $2m\le n$, the space $\H_k$ is
spanned by polynomials of the form $P_k(x)=\det(a'x)^k$ where $a$ is
a complex $n\times m$ matrix satisfying $a'a=0$.
\end{remark}

\noindent {\bf Proof of Theorem \ref{main}.} To prove the first
statement, we consider the equality \be \label {fr1}((T^{\lv}
f)(vr^{1/2}),\; \F\phi)=A(\lv) \; (r_*^{-\lv_*-\nv} \, f(v), \;
\phi),\ee $$ A(\lv)=\frac{(2\pi)^{nm}\,
\Gam_\Om(\lv+\mn)}{c_{\lv}\,\Gam_\Om(-\lv_*)}, $$ which follows from
(\ref{eq9}) and (\ref{fil}). Suppose that $(T^{\lv} f)(v)=0$ almost
everywhere on $\vnm$ for some $\lv \in \frL$. Then $(T^{\lv}
f)(y),\; y=vr^{1/2} \in \Ma$, is zero for almost all $y \in \Ma$,
and (\ref{fr1}) yields $A(\lv)\; (r_*^{-\lv_*-\nv} \, f(v), \;
\phi)=0$. The assumption (\ref{uhh}) together with $\lv \in \frL$
imply $$Re \, \lam_j
>j-m-1 \quad \text{\rm and} \quad \lam_j \neq j-m, j-m+2,\dots \, .$$
Hence $\lv$ is not a pole of $\Gam_\Om(-\lv_*)$, and therefore
$A(\lv) \neq 0$. This gives \be\label {for2}(r_*^{-\lv_*-\nv} \,
f(v), \; \phi)=0,\ee
 the left-hand side being understood in the sense of analytic
 continuation. By choosing $\phi$ as in the proof of Corollary
 \ref{cr}, we obtain $f(v) =0$ a.e. on $\vnm$.

Let us prove the second statement of the theorem. We first note that
for $2m<n$, $H$-polynomials of every degree $k$ do exist. Hence we
can proceed as follows. We observe that the function
$\Gam_\Om(\kv-\lv_*)$ in the expression of $\mu_k(\lv)$
 has the form
\bea \Gam_\Om(\kv-\lv_*)&=&\pi^{m(m-1)/4}\prod\limits_{j=1}^{m}
 \Gam
\Big(\frac{k-(\lam_*)_j- j+1}{2}\Big)\nonumber \\
&=&\pi^{m(m-1)/4}\prod\limits_{j=1}^{m} \Gam \Big(\frac{k+j-\lam_j-
m}{2}\Big)\nonumber \eea (we remind that $(\lam_\ast)_j=Re
\,\lam_{m-j+1}$). It has no poles in $\frL$ if  \be\label{nop}
k>\lam_j +m-j \quad \text{\rm for all} \quad j=1,2,\dots, m.\ee
 Since $\Gam_\Om(\lv+\mn)$ also has no poles in $\frL$, then, by
(\ref{lmu}), $T^\lv P_k =0$ for all $k$ satisfying (\ref{nop}),
provided that $\lv$ is a pole of $\Gam_\Om(-\lv_*)$ (i.e., the
condition (\ref{uhh}) fails). Thus $T^\lv$ is non-injective in this
case. $ \hfill\hfill\square$

We do not know if the condition (\ref{uhh}) is necessary for
injectivity of $T^\lv$ in the case $2m>n$. We suspect that to
answer this question one has to explore $T^\lv$ on polynomial
representations of the  group $SO(n)$ parameterized by highest
weights $(m_1, m_2, \dots , m_{[n/2]})$,  more general  than just
$(k,k,\dots,k)$ corresponding to $\H_k$-spaces; cf. \cite{Str},
\cite{TT}. However, if $\lam_1=\dots = \lam_m=\lam$ then for the
$\lam$-cosine transform
$$(T^{\lam}
f)(u)=\intl_{\vnm} f(v) |\det(v' u)|^{\lam} \, dv\qquad  \Big (
\text{\rm or} \quad (T^{\lam} f)(\xi)=\intl_{G_{n,m}} f(\eta) \,
[\eta | \xi]^{\lam} \, d\eta \Big ),$$ we can give a complete
answer.

\begin{theorem}\label {main2} Let $n> m$ and $Re \,\lam >-1$. Then $T^\lam$
 is injective on
$L^\flat (\vnm) \;$ {\rm(or on $L^1(\gnm)$)} if and only if $\lam
\neq 0,1,2, \dots \,$.
\end{theorem}
\begin{proof} If $\lam_1=\dots \lam_m=\lam$ then the condition
$\lam \neq 0,1,2, \dots$ coincides with (\ref{uhh}). Hence, for
$2m\le n$, the result follows from
 Theorem \ref{main}. If $2m> n$, i.e., $2(n-m)<n$,  we   reduce
 the problem to the dual one on the manifold $V_{n, n-m}$. To this
 end we
  replace the integral $T^\lam f$ by the
equivalent expression in terms of $\gnm$ and write the latter as an
integral over the dual Grassmann manifold $G_{n,n-m}$. This argument
can be realized as follows.

Let $\xi\in\gnm$ and $\eta\in \gnm$ be  $m$-dimensional subspaces
of $\bbr^n$ spanned $u\in \vnm$ and $v \in \vnm$, respectively. We
denote by $\xi^\perp\in G_{n,n-m}$ and $\eta^\perp\in G_{n,n-m}$
the corresponding orthogonal complements, and choose  coordinate
frames $ u_\perp \in V_{n, n-m}$ in $\xi^\perp$ and $ v_\perp \in
V_{n, n-m}$ in $\eta^\perp$. By using the same notation as in
(\ref{ccgr}), we have \be\label{dog}[\eta | \xi]=[\eta ^\perp |
\xi^\perp].\ee Indeed (we recall that $|\cdot|$ denotes the
determinant of the matrix inside)
 \bea
[\eta | \xi]^2&=&|\det(v'u)|^2= |u'vv'u|=|u'(I_n-v_\perp v_{\perp}'
)u|=|I_n-u'v_\perp v_{\perp}' u|\nonumber \\&=&|I_{n-m}- v_{\perp}'
uu'v_\perp|=|I_{n-m}- v_{\perp}'(I_n-u_\perp
u_{\perp}')v_\perp|\nonumber \\&=&| v_{\perp}'u_\perp
u_{\perp}'v_\perp|=|\det( v_{\perp}'u_\perp)|^2\nonumber \\&=&[\eta
^\perp | \xi^\perp]^2.\nonumber\eea Now we successively define the
functions $F(\eta)$ on $\gnm$, $F_\perp(\eta^\perp)$ on $G_{n,n-m}$,
and $f_\perp(v_\perp)$ on $V_{n, n-m}$ by
$$ F(\eta)=f(u), \qquad F_\perp(\eta^\perp)=F((\eta^\perp)^\perp),
\qquad f_\perp(v_\perp)=F_\perp(\eta^\perp).$$ Then
 \bea(T^\lam f)(u)&=&\intl_{G_{n,m}} F(\eta) \, [\eta |
\xi]^{\lam} \, d\eta\nonumber
\\&=&\intl_{G_{n,n-m}}F_\perp(\eta^\perp)\,[\eta
^\perp | \xi^\perp]^\lam \, d\eta^\perp\nonumber \\
&=&\intl_{V_{n,n-m}} f_\perp(v_\perp) |\det(v_{\perp}'
u_\perp)|^{\lam} \,
dv_\perp\nonumber \\
&=&(T_{\perp}^\lam f_\perp)(u_\perp).\nonumber\eea Thus  $T^\lam $
can be expressed as  the similar operator $T_{\perp}^\lam$ on the
dual manifold $V_{n, n-m}$. By above, the latter is non-injective
for $2(n-m)<n$, and we are done. \end{proof}

\bibliographystyle{amsalpha}

\end{document}